\documentclass[11pt]{article}

\usepackage{array}
\usepackage{amsmath}
\usepackage{amssymb}
\usepackage{pstricks}
\usepackage{algorithm}
\usepackage{hyperref}
\usepackage[noend]{algorithmic}

\newcommand{\R}{\mathbb{R}}
\newcommand{\N}{\mathbb{N}}

\newcommand{\C}{\mathbb{C}}

\newcommand{\ac}{\`{a} }

\newtheorem{proposition}{Proposition}
\newtheorem{theorem}{Theorem}
\newtheorem{definition}{Definition}

\newtheorem{lemma}{Lemma}

\def\qed{~\hfill$\Box$\medbreak}

\def\fdimo{~\hfill$\Box$\medbreak}

\begin{document}
	
\title{Large deviation properties for pattern \\ statistics in primitive rational models}

\author{Massimiliano Goldwurm, Marco Vignati}

\date{ }
\maketitle

\begin{center}
{\small Dipartimento di Matematica ``Federigo Enriques'' \\
Universit\ac degli Studi di Milano \\
via Saldini 50, 20133 Milano, Italy \\
	 \{massimiliano.goldwurm,marco.vignati\}@unimi.it}
\end{center}

\begin{abstract}
We present a large deviation property for the pattern statistics representing the number of occurrences of a 
symbol in words of given length generated at random according to a rational stochastic model. 
The result is obtained assuming that in the model the overall weighted transition matrix is primitive.  
In particular we obtain a rate function depending on the main eigenvalue and eigenvectors 
of that matrix.
Under rather mild conditions, we show that the range of validity of our large deviation estimate 
can be extended to the interval (0,1),
which represents in our context the largest possible open interval of validity of the property.
\end{abstract}

\noindent
{\bf Keywords: } regular languages, rational formal series, pattern statistics, large deviations, 
limit distributions.

\section{Introduction}
\label{sec:introd}

Large deviation properties represent a classical subject in probability theory.
They yield bounds of exponential decay 
on the probability that a sequence of random variables 
differs from the mean values for an amount of the order of growth of the mean itself \cite{dz93,dh00}.
Such deviations from the average value are considered ``large'' with respect to other evaluations,
as those deriving for instance from the Central Limit Theorem,
that concern asymptotical smaller differences.

In analytic combinatorics large deviation estimates are  
considered in the study of various relevant structures \cite{fs09}. 
In particular they occur in pattern statistics \cite{js15} 
and in the analysis of depth and height of certain classes of trees \cite{bd06,de99}.
In pattern statistics they have been studied with the goal of evaluating the probability
of rare events, where a given pattern is over- or under-represented in a random text generated according
to a suitable stochastic model \cite{dr04,br14}.

In the present work we prove some properties of this type for sequences of pattern statistics representing
the number of occurrences of a symbol
in a word of length $n$, belonging to a regular language, generated at random according to a
{\em rational stochastic model}.
This model was introduced in \cite{bcgl03} and can be formally 
defined by a nondeterministic finite state automaton
with real positive weights on transitions.
In this setting, the probability of generating a word $w$ of given length is proportional to
the total weight of the accepting transitions labelled by $w$.
This model is quite general, it includes as special cases
the traditional Bernoullian and Markovian sources, widely used in the literature to study the 
number of occurrences of patterns in a random text
\cite{nsf02,rs98,fsv06,br14}. 
We recall that the research concerning pattern statistics has a broad 
range of motivations and applications \cite{js15}.
Moreover, we recall that the rational stochastic model allows also to generate
random words of length $n$ in an arbitrary regular language 
under uniform distribution:
this occurs when the finite automaton defining the model is 
unambiguous and all transitions have the same weight.

In order to fix our notation, consider a (nondeterministic) weighted finite state automaton $\cal A$
over the binary alphabet $\{a,b\}$ and, for every $n\in\N$, let $Y_n$ be
the number of occurrences of the symbol $a$ in a word of length $n$
generated at random according to the rational model defined by $\cal A$.
The analysis of these sequences of random variables is of interest in several contexts.
First of all they can represent the number of occurrences of patterns in a random word of length $n$,
generated by a Markovian source, when the set of patterns is given by a regular language 
\cite{bcgl03,nsf02,rs98}.
Moreover, they are related to the evaluation of the coefficients of rational formal series
(a traditional problem well-studied in the literature \cite{sc62,re77,mr05})
and to the analysis of several problems and properties of regular language.
This fact clearly holds for the natural problem of estimating
the number of words of given length in a regular language having $k$ occurrences of a given symbol
\cite{bmr03,drt00}.
It also holds for the analysis of additive functions defined on regular languages \cite{gr07a}
and for the descriptional complexity of languages and computational models \cite{bmr14}.
Further, using the local limit properties of the sequences $\{Y_n\}$,
for a wide class of rational series it can be proved that the maximum coefficient
of the monomials of degree $n$
has an asymptotic growth of the order $\Theta(n^{k/2}\lambda^n)$ for some $\lambda >0$
and some integer $k \geq -1$ \cite{cgl04,bcgl06}.

The asymptotic behaviour of $\{Y_n\}$, i.e. mean value, variance, 
limit distribution both in the global and in the local sense
\cite{gn97,fs09}, has been studied in the literature under several
hypotheses on the corresponding automaton ${\cal A}$.
It is known that if $\cal A$ has a primitive transition matrix then $Y_n$ has a Gaussian limit distribution 
\cite{bcgl03,nsf02} and, under a suitable aperiodicity condition, 
it also satisfies a local limit theorem \cite{bcgl03},
which can be extended to all primitive cases by using a suitable notion of periodicity \cite{bcgl06}.
The limit distribution of $Y_n$ in the global sense is known also when the transition matrix of $\cal A$ 
consists of two primitive components \cite{dgl04}, 
while the local limit properties in this case are recently studied in \cite{glv21}.
When the automaton ${\cal A}$ has several strongly connected components
a general analysis of the (global) limit distribution of $Y_n$  can be found in \cite{gl06}.

Here we continue this line of research proving in Section \ref{sec:ladeprimo} that
if the transition matrix of ${\cal A}$ is primitive, then
$Y_n$ satisfies a large deviation property with a rate function depending on
the main eigenvalue and the associated eigenvectors.
The corresponding proof is rather standard, it relies on traditional tools of analytic 
combinatorics
and the result is implicitly included in the previous literature \cite{hw96,dh00,fs09}.
However, here our result is significant since it puts in evidence the role played by the main eigenvalue and eigenvectors
of the matrix of weights in the definition of both the rate function and the interval of validity of the property.
Moreover, in Section \ref{sec:interval01},
we assume a mild condition on the transition matrix of the automaton
and, under such hypothesis, we show that the interval of validity of the large deviation property 
can be extended to the entire interval $(0,1)$,
which represents in our context the largest possible open interval of validity of the property.

\section{A quick overview on large deviations}
\label{sec:overview}

Large deviation estimates usually refer to a sequence of random variables, 
say $\{X_n\}$, having increasing mean values;
it consist of a bound, exponentially decreasing to $0$,
over the probability that $X_n$ deviates from $E(X_n)$ by an amount greater or equal to 
$c E(X_n)$, $c>0$. 
Typical situations occur when $E(X_n)\sim \beta n$ for a constant $\beta>0$, and 
since this occurs in all our cases, here we start with
the following fomal definition \cite{dh00,fs09}.

\begin{definition}
{\rm Let $\{X_n\}$ be a sequence of random variables such that $E(X_n) \sim \beta n $ for a constant 
$\beta >0$, and let $(x_0,x_1)$ be an interval including $\beta$.
Assume $I(x)$ is a function defined over $(x_0,x_1)$ taking values in $\R$, 
such that $I(x)> 0$ for $x\neq \beta$.
We say that $\{X_n\}$ satisfies a {\em large deviation property} relative to the interval 
$(x_0,x_1)$ with {\em rate function} $I(x)$
if the following limits hold:
\begin{gather}
\nonumber
\lim_{n\rightarrow \infty} \frac{1}{n} \log \mbox{Pr}(X_n \leq xn) = -I(x) 
\qquad \mbox{ for } x_0 < x \leq \beta \\
\nonumber
\lim_{n\rightarrow \infty} \frac{1}{n} \log \mbox{Pr}(X_n \geq xn) = -I(x) 
\qquad \mbox{ for } \beta \leq x < x_1
\end{gather}
The first relation concerns the {\it left tail} while the second one refers to 
the {\it right tail} of the distribution of $X_n$.}
\end{definition}

This property is equivalent to require that
\begin{gather}
\nonumber
\mbox{Pr}(X_n \leq xn) = e^{-I(x) n + o(n)} \qquad \mbox{ for } x_0 < x \leq \beta \\
\nonumber
\mbox{Pr}(X_n \geq xn) = e^{-I(x) n + o(n)} \qquad \mbox{ for } \beta \leq x < x_1
\end{gather}

A classical example of large deviation property concerns 
the sequence of binomial random variables $\{X_{n,p}\}_n$
of parameters $n$ and $p$, where $p\in (0,1)$ is fixed.
In this case, $E(X_{n,p})=pn$ and by the Central Limit Theorem, we know that
$\frac{X_{n,p}-np}{\sqrt{np(1-p)}}$ converges in distribution to a standard 
Gaussian random variable ${\cal N}(0,1)$.
This yields a limit probability concerning ``normal'' deviations  
(i.e. of the order $\sqrt{n}$) from the mean, that is
$$
\lim_{n\rightarrow \infty} \mbox{Pr}\left(|X_{n,p} - np| \geq \varepsilon \sqrt{n}\right) = 
\mbox{Pr}\left(|{\cal N}(0,1)| \geq\frac{\varepsilon}{\sqrt{p(1-p)}}\right)
\quad \forall \varepsilon > 0
$$
Such a property implies the following result for a larger deviation
$$
\mbox{Pr}\left( |X_{n,p}-np| \geq \varepsilon n \right) = o(1) 
\qquad \forall \varepsilon > 0
$$
which can also be obtained by applying the Law of Large Numbers.
The following proposition proves a large deviation property for $\{X_{n,p}\}_n$ that improves the last relation,
showing that the convergence to $0$ is exponential with respect to $n$
and the range of validity coincides with the overall interval $(0,1)$.

\begin{proposition}
Any sequence of binomial random variables $\{X_{n,p}\}_n$ satisfies a large deviation property 
in the interval $(0,1)$ with rate function $B(x)$ given by
$$
B(x) = x\log \frac{x}{p}+(1-x)\log \frac{1-x}{1-p} \qquad \mbox{ for every } x \in (0,1)
$$
\end{proposition}
{\it Proof.}
First consider the left tail and let $0 < x \leq p$.
We have to prove that
$
\lim_{n\rightarrow \infty} \frac{1}{n} \log \mbox{Pr}(X_{n,p} \leq xn)  = - B(x)
$.
To this end, let $M_n(x) = \max\{ \mbox{Pr}(X_{n,p}=i) : i\in\N, 0 \leq i\leq x n\}$.
Then we have
\begin{equation}
\label{primaeqbin}
M_n(x) \leq \mbox{Pr}(X_{n,p}\leq xn) \leq (xn+1) M_n(x)
\end{equation}
Recall that the probability 
$\mbox{Pr}(X_{n,p} = i) $ is increasing for integers $i$ such that $0 \leq i \leq pn$;
hence
$M_n(x) = {n \choose \lfloor xn \rfloor }p^{\lfloor xn \rfloor}(1-p)^{n-\lfloor xn \rfloor}$
for every $x\in (0,p]$.
Thus, a direct application of Stirling's formula leads to
$$
M_n(x) = {\rm exp}\left\{n \left[x \log \frac{p}{x} + (1-x)\log \frac{1-p}{1-x} \right] +
O(\log n) \right\}
$$
which replaced in (\ref{primaeqbin}) proves that
\begin{equation}
\label{secondaeqbin}
 \log \mbox{Pr}(X_{n,p} \leq xn) = - B(x) n + O(\log n)
\end{equation}
A similar reasoning holds for the right tail.
In this case we have
$$
N_n(x) \leq \mbox{Pr}(X_{n,p}\geq xn) \leq (n - nx +1) N_n(x)
\quad \mbox{ for every } x \in [p,1)
$$
where 
$N_n(x) = 
{n \choose \lfloor xn \rfloor }p^{\lfloor xn \rfloor}(1-p)^{n-\lfloor xn \rfloor}
$.
As above, replacing this value in the previous inequalities yields
\begin{equation}
\label{terzaeqbin}
 \log \mbox{Pr}(X_{n,p} \geq xn) = - B(x) n + O(\log n)
\end{equation}
Relations (\ref{secondaeqbin}) and (\ref{terzaeqbin}) prove the required property.
\qed

The rate function $B(x)$ is strictly convex in the interval $(0,1)$,
takes a unique minimal value at $x=p$ where $B(p)=0$, 
while $B(0^+) = log \frac{1}{1-p}$ and $B(1^-) = log \frac{1}{p}$.
Moreover, it grows vertically for $x$ tending to $0^+$ and to $1^-$, that is
$$\lim_{x\rightarrow 0^+} B'(x) = -\infty \ \mbox{ and }\ 
\lim_{x\rightarrow 1^-} B'(x) = +\infty$$

We recall that often the interval of a large deviation property 
can be extended to the entire set $\R$ once we allow 
the rate function $I(x)$ to assume value $+\infty$.
A classical situation of this type is established by Cram{\'{e}}r's Theorem 
(see for instance \cite{dh00,dz93}), 
stating that if $\{X_n\}$ is a sequence of independent and identically distributed random variables, 
with bounded moment generating function (i.e. $\psi(t)=E(e^{tX_1}) < \infty$ for any $t\in\R$),
then the sequence of partial sums $\{S_n\}$, where $S_n=\sum_{i=1}^nX_i$,
satisfies a large deviation property all over $\R$ with rate function
$$R(x) = \sup_{t\in\R}[xt-\log\psi(t)] \qquad \forall x\in\R$$

\section{Symbol statistics for rational models}

In order to define our stochastic model consider 
a \emph{formal series} in the non-commutative variables $a,b$, that is
a function $r:\{a,b\}^* \rightarrow \R_+$, where $\R_+ =[0,+\infty)$
and $\{a,b\}^*$ is the free monoid of all words on the alphabet $\{a,b\}$.
We denote by $(r,w)$ the value of $r$ at a word $w\in \{a,b\}^*$.
Such a series $r$ is said to be \emph{rational}
if for some integer $m>0$ there exists a 
monoid morphism $\mu : \{a,b\}^* \rightarrow \R_+^{m\times m}$ and two (column) arrays
$\xi, \eta \in \R_+^m$, 
such that $(r,w) = \xi'\mu(w) \eta$, for every $w\in \{a,b\}^*$ \cite{br88,ss78}
(\footnote{As usual we denote by $v'$ the transpose of an array $v\in \R^m$, i.e. a row array.}).
In this case, as the morphism $\mu$ is generated by the matrices
$A=\mu(a)$ and $B=\mu(b)$, we say that the 4-tuple
$(\xi,A,B,\eta)$ is a \emph{linear representation} of $r$.
Clearly, such a 4-tuple can be considered as a finite state automaton over the alphabet $\{a,b\}$,
with transitions weighted by positive real values.
Thus $A$ (resp. $B$) represents the matrix of the weights of all transitions 
labelled by $a$ (resp. $b$), while
$\xi$ (resp. $\eta$) is the array of the weights of the initial (resp. final) states.

Throughout this work, denoting by $\{a,b\}^n$ the family of all words of length $n$ in $\{a,b\}^*$,
we assume that the set $\{w\in\{a,b\}^n : (r,w)>0\}$ is non-empty
for every $n\in\N_+$ (so that $\xi\neq 0\neq \eta$), and that
$A$ and $B$ are non-null matrices.
Thus, we can consider the probability measure $\mbox{Pr}$ over the set $\{a,b\}^n$
given by
$$
\mbox{Pr}(w) = \frac{(r,w)}{\sum_{x\in\{a,b\}^n} (r,x)} =
\frac{\xi'\mu(w)\eta}{\xi'(A+B)^n\eta} \qquad 
\ \forall \ w\in \{a,b\}^n
$$
Note that, if $r$ is the characteristic series of a language $L\subseteq \{a,b\}^*$
then $\mbox{Pr}$ is the uniform probability function over the set
$L\cap \{a,b\}^n$.
Also observe that the traditional Markovian models (to generate a word at random in $\{a,b\}^*$) 
occur when $A+B$ is a stochastic matrix, $\xi$ is a stochastic array and $\eta'=(1,1\ldots,1)$.

Then, under the previous hypotheses, we can define the random variable (r.v.) $Y_n=|w|_a$,
where $w$ is a word chosen at random in $\{a,b\}^n$ with probability $\mbox{Pr}(w)$,
and $|w|_a$ is the number of occurrences of $a$ in $w$. 
As $A\neq [0] \neq B$, $Y_n$ is a non-degenerate random variable.
It is clear that, for every $k\in\{0,1,\ldots,n\}$,
$$
p_n(k) := \mbox{Pr}(Y_n =k) = \frac{\sum_{|w|=n,|w|_a=k} (r,w)}{\sum_{w\in \{a,b\}^n} (r,w)}
$$
Since $r$ is rational also the previous probability can be expressed by using
its linear representation.
It turns out that
\begin{equation}
\label{probipsilon}
p_n(k) = \frac{[x^k]\xi' (Ax+B)^n \eta}{\xi' (A+B)^n \eta}\ ,
\qquad \ k\in\{0,1,\ldots,n\}
\end{equation}
where, as usual, $[x^k] g(x)$ denotes the $k$-th coefficient of the Taylor expansion 
of an analytic function $g(x)$ in a neighbourhood of $0$.

For sake of brevity we say that $Y_n$ is \emph{defined} by the linear representation
$(\xi,A,B,\eta)$.
The moment generating function $\Psi_n(z)$ of $Y_n$ can be defined by means of the map $h_n(z)$ given by
$h_n(z) = \xi' (Ae^z + B)^n\eta$, for $z\in\C$.
We have
\begin{equation}
\label{fgenmom}
\Psi_n(z) = \sum_{k=0}^n p_n(k) e^{zk} = \frac{\xi' (Ae^{z} + B)^n\eta}{\xi' (A+B)^n\eta} = \frac{h_n(z)}{h_n(0)}\ , \qquad \ z \in\C
\end{equation}
and hence mean value and variance of $Y_n$ can be evaluated by
$$
\mbox{E}(Y_n) = \frac{h_n'(0)}{h_n(0)} , \ \ 
\mbox{Var}(Y_n) = \frac{h_n''(0)}{h_n(0)} - \left( \frac{h_n'(0)}{h_n(0)} \right)^2
$$

\section{Primitive models}
\label{sec:primitive}

In this section we resume the main properties of $Y_n$ when the matrix $A+B$ is primitive.
Recall that a matrix $M\in \R_+^{m\times m}$ is {\it primitive} if there exists a positive integer $n$
such that $M^n>0$ (i.e. all entries of $M^n$ are strictly positive).
The main properties of these matrices are established by the following well-known theorem
(see for instance \cite[Sec 1.1]{se81}).
\begin{theorem} (Perron-Frobenius)
\label{teo:pf}
If a matrix $T=[t_{ij}]\in \R_+^{m\times m}$ is primitive then it admits a real eigevalue $\lambda>0$ 
such that:

(i) $|\mu| < \lambda$ for any eigenvalue $\mu$ of $T$ different from $\lambda$;

(ii) $\lambda$ can be associated with strictly positive left and right eigenvectors;

(iii) $\lambda$ is a simple root of the characteristic equation of $T$, 
and hence the associated eigenvectors are unique up to constant multiples;

(iv) if a matrix $A=[a_{ij}]\in \R_+^{m\times m}$ satisfies $A\leq T$ 
(i.e. $a_{ij} \leq t_{ij} \, \forall i,j$) and $\alpha$ is an eigenvalue of $A$ then
$|\alpha| \leq \lambda$. Moreover, $|\alpha| = \lambda$ implies $A=T$.
\end{theorem}
Usually $\lambda$ is called the Perron-Frobenius eigenvalue of $T$.

Then, assume $A+B$ is primitive and let $\lambda$ be its Perron-Frobenius eigenvalue.
In this case it is known that the sequence $\{Y_n\}$ has a Gaussian limit distribution \cite{bcgl03}.
Its properties (in particular mean value and variance) can be studied through the function $y=y(z)$ 
{\em implicitly defined} by the equation 
\begin{equation}
\label{funzipsilon}
\mbox{det}(Iy-Ae^z-B)=0
\end{equation}
with initial condition $y(0)=\lambda$.
Clearly $y(z)$ is eigenvalue of $Ae^z+B$ for every $z\in\C$.
Moreover, $y(z)$ is analytic in a neighbourhood of $0$ and $y'(0)\neq 0$ since
$\lambda$ is a simple root of the characteristic polynomial of $A+B$.

In the analysis of the asymptotic properties of $\{Y_n\}$, 
the following results have been obtained in the literature \cite{bcgl03,bcgl06} 
and are useful in our context:
\begin{description}
	\item[1)] $\mbox{E}(Y_n) = \beta n + c + O(\varepsilon^n)$, where $|\varepsilon|<1$, $c\in \R$ 
	and $\beta$ is a constant satisfying $0<\beta<1$ given by $$\beta = \frac{y'(0)}{\lambda}$$
	Moreover $y'(0)=v'Au$, where $v'$ and $u$ are left and right eigenvectors of $A+B$, with respect to $\lambda$, such that $v'u=1$.
	\item[2)] $\mbox{Var}(Y_n) = \gamma n +O(1)$, where $\gamma$ is a positive constant defined by
	$$\gamma = \frac{y''(0)}{\lambda} - \left(\frac{y'(0)}{\lambda}\right)^2$$
	\item[3)] In a neighbourhood of $0$, the function $\Psi_n(z)$ satisfies a ``quasi power'' condition, 
	that is an equation of the type
	\begin{equation}
	\label{quasiplocale}
	\Psi_n(z) = r(z) \left(\frac{y(z)}{\lambda}\right)^n (1+O(\varepsilon^n)) \qquad (|\varepsilon|<1)
	\end{equation}
	where $r(z)$ is analytic in $z=0$ and $r(0)=1$.
	A consequence of this result is that $\frac{Y_n-\beta n}{\sqrt{\gamma n}}$ 
converges in distribution to a Gaussian random variable of mean $0$ and variance $1$.
\end{description}

Some further properties of the moment generating function $\Psi_n(z)$ can be obtained in the
case of real $z$.
First observe that for every $t\in\R$ also the matrix $Ae^t+B$ is primitive:
therefore $y(t)$ is its Perron-Frobenius eigenvalue.
By the properties of primitive matrices we know that $y(t)$ is a positive real function, 
analytic and strictly increasing for all $t\in\R$ (statement (iv) in Theorem \ref{teo:pf}).
Moreover, all the powers of $Ae^t+B$ satisfy a relation of the form
$$
(Ae^t+B)^n = y(t)^n \cdot u_tv_t' \ (1 + O(\varepsilon_t^n)) \qquad (|\varepsilon_t|<1,\ \forall \ t\in\R)
$$
where $v_t'$ and $u_t$ are left and right eigenvectors of $Ae^t+B$ relative to $y(t)$,
normed so that $v_t'u_t =1$ \cite[Th. 1.2]{se81}.
A first consequence is that applying relation (\ref{quasiplocale}) to all real $z$, 
we obtain (for every $t\in\R$)
\begin{equation}
\label{quasipreale}
\Psi_n(t) = \mbox{E}(e^{tY_n})= \frac{\xi' (Ae^{t} + B)^n\eta}{\xi' (A+B)^n\eta}
= r(t) \left(\frac{y(t)}{\lambda}\right)^n (1+O(\varepsilon_t^n))
\end{equation}
where the function $r(t) = \frac{\xi'u_t v_t' \eta}{\xi'u_0 v_0' \eta}$ is analytic in $\R$,
clearly $r(0)=1$ and $|\varepsilon_t|<1$. 

\section{Large deviations for primitive models}
\label{sec:ladeprimo}

Now assume again $A+B$ primitive and consider the random variable 
$Y_n(t)$ defined by the linear representation $(\xi,Ae^t,B,\eta)$, for any $t\in\R$.
Since $Ae^t+B$ is primitive for any $t\in\R$ we can apply the results of the previous section
to all sequences of random variables $\{Y_n(t)\}$.
To this end, reasoning as for relation (\ref{funzipsilon}), for any $t\in\R$ we can consider the function 
$y_t(z)$ implicitely defined by the equation
$${\rm det}(Iy_t - Ae^{t+z}-B) = 0\ , \qquad z \in\C$$
with initial condition $y_t(0) = y(t)$.
Clearly $y_t(z)=y(t+z)$ and hence $y_t(z)$ is analytic in a neighbourhood 
of $0$ (for any $t\in\R$), 
it admits derivatives of any order around $0$ and 
$y_t'(0) = y'(t)$, $y_t''(0)=y''(t)$.

Applying property {\bf 1)} of the previous section to the linear representation $(\xi,Ae^t,B,\eta)$, 
for every $t\in\R$ we obtain 
$\mbox{E}(Y_n(t)) = \beta(t) n + c_t + O(\varepsilon_t^n)$, 
where $c_t\in \R$ and $\varepsilon_t\in (0,1)$ are constant and $\beta(t)$ is a real function given by
$$
\beta(t) = \frac{y_t'(0)}{y_t(0)} = \frac{y'(t)}{y(t)}  \qquad  \forall t \in \R
$$
Clearly $\beta(0)=\beta$. Moreover, by the same property, we have
$ y(t) = v'_t (Ae^t+B) u_t$, $y'(t) = v'_tAe^tu_t$ 
and hence
\begin{equation}
\label{betadit2}
\beta(t) = \frac{v'_tAe^tu_t}{v'_t (Ae^t+B) u_t} 
\end{equation}
with 
\begin{equation}
\label{betadit3}
0 < \beta(t) < 1  \qquad  \forall t \in \R
\end{equation}

Analogously, applying property {\bf 2)} of the previous section to $Y_n(t)$ we get
$$
\mbox{Var}(Y_n(t)) = \gamma(t) n + O(1)  \qquad \forall t \in \R
$$
where $\gamma(t)$ is a positive constant given by
\begin{equation}
\label{gammaditipos}
\gamma(t) =  \frac{y_t''(0)}{y_t(0)} - \left(\frac{y_t'(0)}{y_t(0)}\right)^2 = 
\beta'(t) > 0  \qquad \forall t \in \R
\end{equation}
Therefore $\beta(t)$ is strictly increasing all over $\R$ and
the following limits exist and are finite:
\begin{equation}
\label{limiti}
U = \lim_{t\rightarrow -\infty}\beta(t) ,\qquad  V=\lim_{t\rightarrow +\infty}\beta(t)
\end{equation}
By relation (\ref{betadit3}), we have
$$0 \leq U < \beta(0) < V \leq 1$$
which, together with relation (\ref{gammaditipos}), implies the following statement.

\begin{lemma}
\label{tauconx}
For every $x\in (U,V)$ there exists a unique $\tau_x \in \R$ such that 
\begin{equation}
\label{equaztau}
\beta(\tau_x)=x
\end{equation}
Moreover, $\tau_x < 0$ whenever $x<\beta$, $\tau_{\beta}=0$ and $\tau_x> 0$ when $x>\beta$.
\end{lemma}

Now we apply property {\bf 3)} of the previous section to the random variable
$Y_n(t)$:
we get a ``quasi power'' property for the moment generating function of $Y_n(t)$,
that is
$\Psi_{Y_n(t)}(z) = r_t(z) \left(\frac{y(t+z)}{y(t)}\right)^n (1+O(\varepsilon_t^n))$,
for some $\varepsilon_t\in (0,1)$, where $r_t(z)$ is also analytic in $z=0$ and $r_t(0)=1$.
As a consequence, for every $t\in\R$ the sequence of random variables 
$\left\{\frac{Y_n(t)-\beta(t) n}{\sqrt{\gamma(t) n}}\right\}_n$ 
converges in distribution to a Gaussian random variable of mean $0$ and variance $1$,
i.e. for every constant $x\in\R$ we have
\begin{equation}
\label{convergenzaydit}
\lim_{n\rightarrow \infty} \mbox{Pr}\left( \frac{Y_n(t)-\beta(t) n}{\sqrt{\gamma(t) n}} \leq x \right) 
= \frac{1}{\sqrt{2\pi}} \int_{-\infty}^x e^{-t^2/2} dt \qquad \forall t\in\R
\end{equation}

The previous results allows us to prove a large deviation property for $\{Y_n\}$.
\begin{theorem}
\label{teo:ladevYdin}
	Let $\{Y_n\}$ be defined by a linear representation $(\xi,A,B,\eta)$ where $A+B$ is primitive.
	Then $\{Y_n\}$ satisfies a large deviation property in the interval $(U,V)$ with rate function
	$I(x) = - \log \left( \frac{y(\tau_x)}{\lambda e^{x\tau_x}}\right)$,
	where $\tau_x$ is defined by equation (\ref{equaztau}).
\end{theorem}
{\bf Proof.}
We first study the right tail of $\{Y_n\}$.
We have to prove that for every $x\in [\beta,V)$ the following relation holds:
\begin{equation}
\label{rightail}
\lim_{n\rightarrow +\infty} \frac{1}{n} \log \mbox{Pr}(Y_n \geq xn) = 
\log \left( \frac{y(\tau_x)}{\lambda e^{x\tau_x}}\right)
\end{equation}
By Markov inequality, for every $t>0$ we have
$$
\mbox{Pr}(Y_n \geq xn) = \mbox{Pr}(e^{tY_n} \geq e^{txn}) \leq \frac{\mbox{E}(e^{tY_n})}{e^{txn}}
$$
and hence, by relation (\ref{quasipreale})
we get
$\mbox{Pr}(Y_n \geq xn) \leq r(t) \left(\frac{y(t)}{\lambda e^{tx}}\right)^n (1+O(\varepsilon_t^n))$,
which implies
$$
\frac{1}{n} \log \mbox{Pr}(Y_n \geq xn) \leq
\log \left( \frac{y(t)}{\lambda e^{xt}}\right) +O(1/n)
$$
This bound can be further refined by taking the minimum with respect to $t>0$ 
of the first term in the right hand side.
To this end let us define the function 
\begin{equation}
\label{phiconx}
\varphi_x(t) = \log \left( \frac{y(t)}{\lambda e^{xt}}\right) \qquad \forall t\in \R
\end{equation}
Note that $\varphi_x(0)=0$, $\varphi_x'(t) = \beta(t) - x$, and hence by Lemma \ref{tauconx}
since $x\geq \beta$, $\varphi_x(t)$ takes a unique minimum value at $t=\tau_x \geq 0$.
This proves
$$
\lim_{n\rightarrow +\infty} \frac{1}{n} \log \mbox{Pr}(Y_n \geq xn) \leq
\log \left( \frac{y(\tau_x)}{\lambda e^{x\tau_x}}\right)
$$
Also observe that $\varphi_x(t)$ is a convex function since $\varphi_x''(t)=\beta'(t) >0$ 
by relation (\ref{gammaditipos}).

An analogous lower bound for $\mbox{Pr}(Y_n \geq xn)$ can be proved 
by considering the random variable $Y_n(\tau_x)$.
Since
$[z^k] \xi'(Ae^{\tau_x}z+B)^n\eta = e^{\tau_x k} [z^k] \xi'(Az+B)^n \eta$, 
by relations (\ref{probipsilon}) and (\ref{fgenmom}) we have
\begin{equation}
\label{probycontau}
\mbox{Pr}(Y_n=k) = 
\frac{ \mbox{Pr}\left(\begin{array}{c} Y_n(\tau_x)=k \end{array} \right)\; \Psi_n(\tau_x) }{e^{\tau_x k}} 
\qquad \forall \ k =0,1,\ldots,n
\end{equation}
Also note that $\mbox{E}(Y_n(\tau_x)) = \beta(\tau_x)n + O(1) = xn +O(1)$ and by (\ref{convergenzaydit})
we know that $\{Y_n(\tau_x)\}_n$ has a Gaussian limit distribution.
This means that, for every $\varepsilon >0$,
$\mbox{Pr}\left(\ \begin{array}{c} Y_n(\tau_x) > (x+\varepsilon)n \end{array}\  \right) =o(1)$ and then
\begin{equation}
\mbox{Pr}\left(\  \begin{array}{c} xn \leq Y_n(\tau_x) \leq (x+\varepsilon)n \end{array} \ \right) = 
\frac{1}{2} + o(1)
\end{equation}
Then, from this relation and the identities (\ref{probycontau}) and (\ref{quasipreale}) we get
\begin{eqnarray}
\nonumber
\mbox{Pr}(Y_n \geq xn) & \geq & 
\mbox{Pr}\left(\ \begin{array}{c} xn \leq Y_n \leq (x+\varepsilon)n \end{array} \ \right)  \\
\nonumber
& \geq & 
\frac{\mbox{Pr}\left(\ \begin{array}{c} xn \leq Y_n(\tau_x) \leq (x+\varepsilon)n\end{array} \ \right) \ \ \Psi_n(\tau_x)}{e^{\tau_x(x+\varepsilon)n}} = \\
& = & \left( \frac{1}{2} + o(1) \right) r(\tau_x)
\left(\frac{y(\tau_x)}{\lambda e^{\tau_x(x+\varepsilon)}}\right)^n (1+O(\varepsilon_{\tau_x}^n))
\label{3stars}
\end{eqnarray}
Thus, by the arbitrariness of $\varepsilon$, we have
$$
\frac{1}{n} \log \mbox{Pr}(Y_n \geq xn) \geq \log \left( \frac{y(\tau_x)}{\lambda e^{x\tau_x}}\right) + O(1/n)
$$
which yields the required lower bound and concludes the proof of relation (\ref{rightail}).

Consider now the left tail. We have to prove that,
for every $x\in (U,\beta ]$,
\begin{equation}
\label{leftail}
\lim_{n\rightarrow +\infty} \frac{1}{n} \log \mbox{Pr}(Y_n \leq xn) = 
\log \left( \frac{y(\tau_x)}{\lambda e^{x\tau_x}}\right)
\end{equation}
where $\tau_x \leq 0$ is defined by equation (\ref{equaztau}).

The reasoning is similar to the previous case.
The main difference is that here $U < x \leq \beta$ and one has to use negative values of $t$.
Note that the function $\phi_x(t)$ given by (\ref{phiconx}) is well defined also in this case.
For every $x\in (U,\beta ]$ and every $t<0$, by Markov inequality and relation (\ref{quasipreale})
we get 
$$
\mbox{Pr}(Y_n \leq xn)  = \mbox{Pr}(e^{tY_n} \geq e^{txn}) \leq \frac{\mbox{E}(e^{tY_n})}{e^{txn}} 
=  r(t) \left(\frac{y(t)}{\lambda e^{tx}}\right)^n (1+O(\varepsilon_t^n))
$$
By Lemma \ref{tauconx} the minimum of 
$\varphi_x(t) = \log \left( \frac{y(t)}{\lambda e^{xt}}\right)$ 
is taken at $t=\tau_x \leq 0$ and this 
proves that
$$
\frac{1}{n} \log \mbox{Pr}(Y_n \leq xn) \leq 
\log \left( \frac{y(\tau_x)}{\lambda e^{x\tau_x} } \right) + O(1/n)
$$
which yields an upper bound to the limit in (\ref{leftail}).

The corresponding lower bound is obtained as in the analysis of the right tail,
leading to relation (\ref{3stars}),
with obvious changes.
\fdimo

\section{Large deviations in the interval (0,1)}
\label{sec:interval01}

A natural question arising at this point is whether the interval $(U,V)$ can be extended to $(0,1)$ 
as in the case of the sequences of binomial random variables considered in section \ref{sec:overview}.

Let us assume that $A+B$ is a primitive matrix. Since $A$ and $B$ are non-null matrices
with entries in $\R_+$, they admit a real non-negative eigenvalue that is greater or equal to the 
modulus of any other eigenvalue of the respective matrix.
We denote such eigenvalues by $\lambda_A$ and $\lambda_B$, respectively.
Clearly, as $A$ is not primitive in general, it may occur $\lambda_A=0$ or $\lambda_A=|\mu|$ 
for some eigenvalue $\mu$ of $A$ such that $\mu\neq \lambda_A$,
and the same may happen for $\lambda_B$.
However, by statement (iv) of Theorem \ref{teo:pf}, it is clear that
$\lambda_A < \lambda$ and $\lambda_B < \lambda$, where $\lambda$ is the Perron-Frobenius
eigenvalue of $A+B$.

Now, assume $\lambda_B>0$ (which is equivalent to require that $B$ has a nonnull eigenvalue)
and let $v_B$ and $u_B$ be left and right eigenvectors of $B$
with respect to $\lambda_B$, normed so that $v'_B u_B = 1$.
Clearly $v_B$ and $u_B$ cannot be null.
Moreover, for $t\rightarrow -\infty$, the matrix $Ae^t + B$ tends to $B$
and hence the eigenvalue $y(t)$ converges to $\lambda_B$, while the matrix $u_t v'_t$ tends 
$u_B v'_B$, implying $v'_tBu_t \rightarrow v'_B B u_B=\lambda_B$ and similarly
$v'_t A u_t \rightarrow v'_B A u_B$.
As a consequence, for $t\rightarrow -\infty$, we have
$$
y'(t) = v'_tAe^tu_t = O(e^t) = o(1)
$$
Therefore, by equality (\ref{betadit2}) the last relation implies (for $t\rightarrow -\infty$)
\begin{equation}
\label{limiteameno}
\beta(t) = \frac{v'_tAe^tu_t}{v'_t(Ae^t+B)u_t} = \frac{O(e^{t})}{\lambda_B + o(1)} = o(1)
\end{equation}
and hence $U=0$.

Analogously, if $\lambda_A>0$ we get $V=1$.
In fact, assume $t\rightarrow +\infty$ and let $L(t)= A+ Be^{-t}$.
Exchanging $A$ and $B$ in the previous argument and recalling that $L(t)$ and 
$Ae^t + B$ have the same eigenvectors, we obtain
$v'_t A u_t = \lambda_A + o(1)$ and $v'_tBe^{-t}u_t = O(e^{-t}) = o(1)$.
As a consequence, for $t\rightarrow +\infty$, we have
\begin{equation}
\label{limiteapi}
\beta(t) = \frac{v'_tAe^tu_t}{v'_t(Ae^t+B)u_t} = 
\frac{v'_tAu_t}{v'_t(A+Be^{-t})u_t} = 
\frac{\lambda_A + o(1)}{\lambda_A + o(1)} = 1 +o(1)
\end{equation}
which implies $V=1$.

The previous argument proves the following result.
\begin{theorem}
	\label{teo:ladin01}
	Let $\{Y_n\}$ be defined by a linear representation $(\xi,A,B,\eta)$ where 
	$A+B$ is primitive and both $A$ and $B$ have a nonnull eigenvalue.
	Then $\{Y_n\}$ satisfies a large deviation property in 
	the interval $(0,1)$ with rate function
	$$I(x) = - \log \left( \frac{y(\tau_x)}{\lambda e^{x\tau_x}}\right)\qquad \forall x\in (0,1)$$
	where $\tau_x$ is the unique real value such that $\beta(\tau_x) = x$.
\end{theorem}

Under the same hypotheses of the previous theorem we can study the function $I(x)$.
Note that the function $\tau=\tau_x$, implicitly defined
by $\beta(\tau)-x = 0$, is well defined and analytic for $x\in (0,1)$.
Thus it is easy to see that
$$
I'(x) = - \beta(\tau_x) \tau'_x + \tau_x + x\tau'_x = \tau_x
$$
and hence $I(x)$ is decreasing in $(0,\beta)$ and increasing in $(\beta,1)$,
with a unique minimal value $I(\beta)=0$.
This also proves that 
$$\lim_{x\rightarrow 0^+} I'(x) = -\infty \mbox{ and }
\lim_{x\rightarrow 1^-} I'(x) = +\infty$$
and hence $I(x)$ grows vertically for $x\rightarrow 1^-$ and for $x \rightarrow 0^+$.
Moreover, $I''(x)$ equals $\tau'_x$, which is always positive since $\tau_x$ is strictly increasing
in $(0,1)$, and hence $I(x)$ is a convex function in $(0,1)$.

Finally, let us determine the behaviour $I(x)$ for $x \rightarrow 0^+$ and for
$x \rightarrow 1^-$.
Letting $x \rightarrow 0^+$, we have $\tau_x \rightarrow -\infty$ and,
arguing as for equation (\ref{limiteameno}), we obtain
$y(\tau_x) = \lambda_B + o(1)$.
Moreover, applying the same equation (\ref{limiteameno}), we get
$$x\tau_x = \beta(\tau_x) \tau_x = O(e^{\tau_x}) \tau_x = o(1)$$
As a consequence, we can derive the following limit
$$
\lim_{x \rightarrow 0^+} I(x) = 
\lim_{x \rightarrow 0^+} -\log(y(\tau_x)) + \log(\lambda)+ x\tau_x = 
\log (\lambda/\lambda_B)
$$

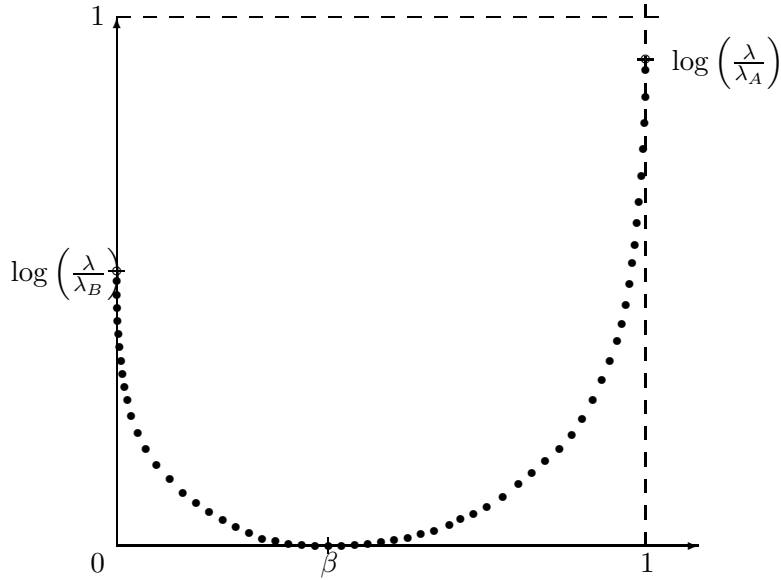
\begin{figure}[h]
	\begin{center}
		\begin{picture}(300,220)(0,0)
			\put (40,20){\vector(1,0){220}}
			\put (30,10){$0$}
			\put (238,10){$1$}
			\multiput (240,20)(0,10){21}{\line(0,1){5}}	
			\put (117,10){$\beta$}
			\put (120,17){\line(0,1){6}}
			\put (40,20){\vector(0,1){200}}
			\put (30,217){$1$}
			\multiput (40,220)(10,0){21}{\line(1,0){5}}
			
			\put(250,200){$\log \left(\frac{\lambda}{\lambda_A}\right)$}
			\put (243,204){\line(-1,0){6}}
			\put(0,120){$\log \left(\frac{\lambda}{\lambda_B}\right)$}
			\put (37,124){\line(1,0){6}}
			
			\put (40,124){\circle{3}}
			
			\put (40,120){\circle*{3}}
			\put (40.1,115){\circle*{3}}
			\put (40.2,110){\circle*{3}}
			\put (40.3,105){\circle*{3}}
			\put (40.6,100){\circle*{3}}
			\put (41,95){\circle*{3}}
			\put (41.6,90){\circle*{3}}
			\put (42.1,85){\circle*{3}}
			\put (43,80){\circle*{3}}
			\put (44,75){\circle*{3}}
			
			\put (45.5,69){\circle*{3}}
			\put (48,62.6){\circle*{3}}
			\put (51,56.65){\circle*{3}}
			\put (55,50.5){\circle*{3}}
			
			\put (60,45){\circle*{3}}
			\put (65,40){\circle*{3}}
			\put (70,36.1){\circle*{3}}
			\put (75,32.5){\circle*{3}}
			\put (80,29.5){\circle*{3}}
			\put (85,27){\circle*{3}}
			\put (90,24.6){\circle*{3}}
			\put (95,22.7){\circle*{3}}
			\put (100,21.7){\circle*{3}}
			\put (105,20.65){\circle*{3}}
			\put (110,20.40){\circle*{3}}
			\put (115,20.05){\circle*{3}}
			\put (120,20){\circle*{3}}
			
			\put (125,20){\circle*{3}}
			\put (130,20.3){\circle*{3}}
			\put (135,20.6){\circle*{3}}
			\put (140,21.3){\circle*{3}}
			\put (145,22){\circle*{3}}
			\put (150,23){\circle*{3}}
			\put (155,24.2){\circle*{3}}
			\put (160,25.5){\circle*{3}}
			\put (165.8,27.9){\circle*{3}}
			\put (170,30){\circle*{3}}
			
			\put (175,32){\circle*{3}}
			\put (180,34.7){\circle*{3}}
			\put (186,38.4){\circle*{3}}
			\put (191.9,43.1){\circle*{3}}
			\put (197,47.3){\circle*{3}}
			\put (202,52){\circle*{3}}
			\put (207.5,56.5){\circle*{3}}
			\put (212,62){\circle*{3}}
			
			\put (216,68){\circle*{3}}
			\put (220,75){\circle*{3}}
			\put (223.5,82.5){\circle*{3}}
			\put (226.5,90){\circle*{3}}
			\put (229.3,97.5){\circle*{3}}
			\put (231.1,104){\circle*{3}}
			\put (232.5,111){\circle*{3}}
			\put (234,119){\circle*{3}}
			
			\put (234.9,127){\circle*{3}}
			
			\put (235.9,134){\circle*{3}}
			\put (236.7,142){\circle*{3}}
			
			\put (237.5,150){\circle*{3}}
			\put (238.4,160){\circle*{3}}
			\put (239.2,170){\circle*{3}}
			\put (239.6,180){\circle*{3}}
			\put (240,190){\circle*{3}}
			\put (240,200){\circle*{3}}
			
			\put (240,204){\circle{3}}
		\end{picture}
		\caption{Function $I(x)$ with $\beta =\frac{v'Au}{\lambda} = 2/5$, assuming $\lambda_A=v'Au$
		and $\lambda_B=v'Bu$.}
		\label{ratefunc}
	\end{center}
\end{figure}

Analogously, let $x \rightarrow 1^-$.
Then $\tau_x \rightarrow +\infty$ and, reasoning as for equation (\ref{limiteapi}), we get 
$y(\tau_x) = e^{\tau_x}( \lambda_A + o(1))$ and
$\beta(\tau_x)=1+O(e^{-\tau_x)}$.
This implies 
$$x\tau_x = \beta(\tau_x) \tau_x = \tau_x + o(1) $$
and hence we obtain the following limit
$$
\lim_{x \rightarrow 1^-} I(x) = 
\lim_{x \rightarrow 1^-} -\log(y(\tau_x)) + \log(\lambda)+ x\tau_x = 
\log (\lambda/\lambda_A)
$$
The diagram of an example of function $I(x)$ is described in Figure \ref{ratefunc}.

In conclusion we can see that the rate function $I(x)$ is rather 
similar to the map $B(x)$ discussed in section \ref{sec:overview}, 
that is the rate function of the large deviation property 
for the sequence $\{X_{n,p}\}$ of binomial random variables.
In particular the constant $\beta=\frac{v'Au}{\lambda}$, where $I(x)$ takes its minimal
value $0$, seems to play the same role as the success probability $p$ of $X_{n,p}$.
Also observe that the limits of $I(x)$ for $x\rightarrow 0^+$
and $x\rightarrow 1^-$, i.e. $\log (\lambda/\lambda_B)$ and $\log (\lambda/\lambda_A)$ respectively,
correspond to $-\log (1-p)$ and $-\log p$,
which seems to suggest that $\lambda_A= v'Au$ and $\lambda_B=v'Bu$.
%

\newpage
\bibliographystyle{plain}

\end{document}